\documentclass[12pt, a4paper]{amsart}
\usepackage{amssymb}
\usepackage{dsfont}
\newtheorem{proposition}{Proposition}[section]
\newtheorem{theorem}{Theorem}[section]

\newtheorem{corollary}{Corollary}[section]
\newtheorem{definition}{Definition}[section]
\newtheorem{remark}{Remark}[section]

\newcommand{\eps}{\varepsilon}
\newcommand{\R}{{\mathbb{R}}}
\newcommand{\N}{{\mathbb{N}}}
\renewcommand{\P}{{\mathds{P}}}

\begin{document}

\newcommand{\dmu}{\,\mathrm{d}\mu}

\title{The Dirichlet problem for some nonlocal diffusion equations}

\date{}
\author{Emmanuel Chasseigne${}^1$}

\footnotetext[1]{Laboratoire de Math\'ematiques et Physique Th\'eorique, Universit\'e F. Rabelais, Parc de Grandmont, 37200 Tours (FRANCE) - e-mail : emmanuel.chasseigne@lmpt.univ-tours.fr}

\begin{abstract}
We study the Dirichlet problem for the non-local diffusion equation
$u_t=\int\{u(x+z,t)-u(x,t)\}\dmu(z)$, where $\mu$ is a $L^1$ function and $``u=\varphi$ on $\partial\Omega\times(0,\infty)$'' has to be understood in a non-classical sense. We prove existence and uniqueness results of solutions in this setting. Moreover, we prove that our solutions coincide with those obtained through the standard ``vanishing viscosity method'', but show that a boundary layer occurs: the solution does not take the boundary data in the classical sense on $\partial\Omega$, a phenomenon related to the non-local character of the equation. Finally, we show that in a bounded domain, some regularization may occur, contrary to what happens in the whole space.
\end{abstract}
\maketitle

\date{}

\medskip


\section{Introduction}
\setcounter{equation}{0}

In this paper we investigate a suitable notion of Dirichlet Problems in a bounded set $\Omega\subset\R^n$ for a class of non-local equations, starting from remarks and ideas developed in \cite{CCR} to cover a more general framework. 

To begin with, let us consider our model equation posed in the whole space $\R^n$ :
\begin{equation}\label{eq0}\frac{\partial u}{\partial t}(x,t)=\int_{\R^n}\big\{u(x+z,t)-u(x,t)\big\}\dmu(z)\,,
\end{equation}
where $\mu\in\mathcal{M}(\R^n)$ is a finite measure. Such non-local diffusion operators arise in various problems like diffusion of population (see \cite{CCR} and the references therein for more on this topic), L\'evy processes \cite{bertoin,CV,woyczynski}, fractional Laplace diffusion \cite{silvestre06}, and more generally, integro-differential pde's (see \cite{Al} and the references therein). Note that in the fractional Laplacian case,
$\mu(z)=|z|^{-N-\lambda}$ with $0<\lambda<2$, a case in which $\mu$ is not a finite measure that we shall not consider here (we refer to \cite{BCI} for a general study of the Dirichlet problem in such situations). 

Now, if the equation has a rather clear sense when posed in the entire space $\R^N$, the question of what is a Dirichlet problem in this context becomes less obvious. Indeed, due to the nonlocal character of the equation, we need to assign the values of $u$ not only on the topological boundary $\partial\Omega$, but on the complement of $\Omega$ (see more precise statments below). 

Let us mention that a viscosity approach of such Dirichlet problems can be found in \cite{Ar} and \cite{CV}, but it concerns mainly the case of singular (\textit{i.e.}, non-integrable) measures. We shall prefer here a similar approach to that which can be found in \cite{CCR} where we considered an homogenous Dirichlet problem with $\mu$, a continuous bounded symetric function. But we want to go further in two directions: first to consider non-homogeneous problems (to be suitably defined), and second to include functions $\mu$ that may be singular, but yet integrable.

\

\noindent\textsc{Aims and Organization of the Paper}

The aim of this paper is thus to give a rather complete view of the case when $\mu\ge0$ is an integrable function in $\R^n$. 

We shall first define a suitable notion of what is a Dirichlet problem for our equation in Section 2, and then we investigate solvability of the problem in Section 3. More precisely, we prove that for any continuous and bounded initial data $u_0$ and ``boundary'' data $\varphi$, there exists a unique solution $u$, continuous and bounded in $\bar\Omega\times[0,\infty)$. We also prove that comparison holds in this class of solutions.

Section 4 is devoted to study the vanishing viscosity method which is another way to try to solve the problem :  we add a viscous term $-\epsilon\Delta u$ to the equation in order to force the boundary data to be taken, and then try to see what happens in the limit as $\epsilon\to0$. Well, it turns out that we get indeed in the limit the unique solution obtained in Section 3, but this solution is strictly positive on $\partial\Omega$ whereas $u_\epsilon$ is zero on the boundary for any $\epsilon>0$. Hence, we are facing a boundary layer formation which proves that in general the boundary data will not be taken in the usual sense.

In Section 5 we investigate some fine regularity properties of the solutions. It is known that in the case $\Omega=\R^n$, no regularization occurs for such non-local equations, in the sense that the solution at time $t>0$ cannot have a  better regularity than that of $u_0$ (this was proven in \cite{CCR}). But in the case of bounded domains, it appears that however some kind of smoothing does occur, provided $\mu$ is not too much concentrated: the modulus of continuity of the solution decreases exponentially as time increases (see Section 5 for precise statments).

We will focus on some other qualitative aspects in Section 6. In particular, we prove that the equation enjoys an positivity property similar to the heat equation (if $u_0\ge0$ and not identically zero, then $u(t)$ is automatically positive in all $\Omega$). But in the present situation, this effect is not purely diffusive as for the Laplace operator, it is due to the nonlocal character of the equation. Moreover we show that solutions are immediately positive also on the boundary $\partial\Omega$, which is consistent with the layer formation in the vanishing viscosity method. 

Finally we collect some other qualitative remarks, examples and possible extensions in the last Section.

\

\noindent\textit{This note benefited a lot from discussions with Guy Barles and Ch. Georgelin, whom I warmly thank.}

\section{Definition of the Dirichlet Problem}

We shall assume throughout this paper that the measure $\mu$ has a density (still denoted by $\mu$) which is a nonnegative, integrable fonction of total mass $\mu(\R^n)=1$. In fact, if the density is integrable, we can always assume that the total mass is $1$ by invoking scaling arguments: simply put $v(x,t)=u(x,\lambda t)$ where $\lambda$ is the total mass of $\mu$. The set $\Omega\subset\R^n$ will be open, bounded, connected and regular (to simplify).

Now, as we noticed above, to give sense to a Dirichlet problem one needs to prescribe the ``boundary'' value not only on the topological frontier $\partial\Omega$, but in fact in the complement of $\Omega$ or to be more precise, on the set where the measure $\mu$ will read some information. So let $\tau^{-1}_x$ be the translated measure defined by $\tau^{-1}_x\mu(z)=\mu(x+z)$, let $\mathrm{supp}\,(\mu)$ denote the support of the measure $\mu$ and let us define the following subsets of $\R^n$ :

$$\bar\Omega^\mu = \overline{\bigcup_{x\in\Omega}\mathrm{supp}\,(\tau^{-1}_x\mu)}\,,\ \mbox{and}\ \partial^{\mu}\Omega=\bar\Omega^{\mu}\setminus\Omega\,.$$ 

In other words, $\bar\Omega^{\mu}$ is what we call the \textit{domain of influence} of $\Omega$. In short, if you want to know how things evolve in $\Omega$, you need to read information in the bigger set $\bar\Omega^\mu\supset\Omega$, because the operator will look at it. In this direction, the natural notion of frontier $\partial\Omega$ has to be replaced by the \textit{extended frontier} $\partial^\mu\Omega$, which is of course not a negligeable set since it is regular (\textit{i.e.}, the closure of an open set in $\R^N$). Then in this setting, we may construct a general theory for the following Dirichlet problem :  
$$(\mathcal{P})\quad\begin{cases}
   \dfrac{\partial u}{\partial t}(x,t)=\displaystyle\int_{\R^n}\big\{u(x+z,t)-u(x,t)\big\}\dmu(z) & \mbox{in}\ \Omega\times(0,\infty)\,, \cr
   u(x,t)=\varphi(x,t) & \mbox{in}\ \partial^\mu\Omega\times[0,\infty)\,,\cr
   u(x,0)=u_0(x) & \mbox{in}\ \Omega\,.
\end{cases}
$$

The preceeding formulation requires the function $u$ to be defined in $\Omega^\mu\times[0,\infty)$, extended by $\varphi$, in order to give sense to the integral term. But as we shall see below, $u$ will not be continuous when crossing the boundary $\partial\Omega$ so that practically, we shall often prefer a different formulation, keeping the function $u$ to be defined only within $\Omega\times[0,\infty)$. Let us denote by $K_\varphi(u)$ the following non-local operator : for all $(x,t)\in\Omega\times(0,\infty),$ 
$$K_\varphi(u)(x,t)=\int_{\{(x+z)\in\Omega\}}u(x+z,t)\dmu(z)+\int_{\{(x+z)\in\partial^\mu\Omega\}}\varphi(x+z,t)\dmu(z)-u(x,t)\,.$$

Thus we arrive at the following definition :

\begin{definition}\label{def}
Let $\mu$ be a finite measure, $\varphi\in\mathrm{C}(\partial^\mu\Omega\times[0,\infty))$ and $u_0\in \mathrm{C}(\bar{\Omega})$ with $u_0=\varphi$ in $\partial\Omega\times\{0\}$. By a solution of $(\mathcal{P})$, we mean a function $u\in\mathrm{C}(\Omega\times[0,\infty[)$ such that 
$u_t=K_\varphi(u)$ in $\Omega\times(0,\infty)$ in the sense of distributions and $u(0)=u_0$ in $\Omega$.
\end{definition}

Note that writing $u_t=K_\varphi(u)$ implicitly means that we are trying to solve the Dirichlet problem with ``$u=\varphi$'' in $\partial^\mu\Omega\times(0,\infty)$, while $u$ is only defined in $\Omega\times[0,\infty)$. Note also that since we consider continuous solutions, $u_t$ will be continuous and the equation will hold in the classical sense.

A more general definition could be given, only imposing that the translations of $\varphi$ and $u_0$ be integrable with respect to $\mu$. Thus, in the case when $\mu$ is a bounded function, it is enough to consider that $\varphi$ and $u$ are integrable with respect to the Lebesgue measure and the formulation makes sense. But in this paper, we shall restrict ourselves to the case $\mu\in L^1(\R^n)$ (see next Section). 

Notice that if the support of $\mu$ is not compact, then we may have possibly $\bar\Omega^\mu=\R^N$ and $\partial^\mu\Omega=\R^N\setminus\Omega$, the case in which every point, even outside $\Omega$ will influence the density at any other point $x$ in $\Omega$. On the contrary, if $\mu$ is compactly supported, then the values of $u$ outside $\bar\Omega^\mu$ will never enter into play, as soon as we look at the problem in $\Omega$. Of course the homogeneous case treated in \cite{CCR} appears as a special case of this framework when $\varphi\equiv0$ and $\mu\in\mathrm{C}(\R^n)$. 

We may now come back to our first intent and make precise the formulation : if $u$ is a solution in the sense of Definition \ref{def}, we can extend it as follows
\begin{equation}
\tilde{u}(x,t)=\begin{cases}\label{extension}
u(x,t) & \text{if }x\in\Omega\,,\\ \varphi(x,t) & \text{if }x\in\partial^\mu\Omega\,,
\end{cases} 
\end{equation}
and for instance $\tilde{u}=0$ outside $\Omega^\mu\times(0,\infty)$ (the values there are not involved in the integral). Then $\tilde{u}$ satisfies :
$$\frac{\partial \tilde{u}}{\partial t}(x,t)=\displaystyle\int_{\R^n}\big\{\tilde{u}(x+z,t)-\tilde{u}(x,t)\big\}\dmu(z)\ \mbox{in}\ \Omega\times(0,\infty)\,.$$

\section{A General Theory of Existence and Uniqueness}

In this Section, we shall derive here a general existence and uniqueness result when $\mu$ is a $L^1$ measure. In the following, the notation $\mathrm{C}_b(Y)$ means the set of continuous and bounded functions on $Y$, and let $$X=\mathrm{C}_b(\bar\Omega\times[0,\infty))\,.$$

\begin{theorem}\label{thmexun}
Let $\mu\in L^1(\R)$. Then for any $u_0\in\mathrm{C}_b(\bar\Omega)$ and $\varphi\in \mathrm{C}_b(\partial^\mu\Omega\times[0,\infty))$ with $u_0=\varphi$ in $\partial\Omega\times\{0\}$, there exists a unique solution $u\in X$ of the Dirichlet problem $(\mathcal{P})$. 
\end{theorem}

\begin{proof}
This result is obtained by a standard Banach fixed-point theorem in $X$ as follows: consider for $t_0>0$ fixed the Banach space $X_{t_0}=\mathrm{C}([0,t_0]\times\bar\Omega)$ and the operator $T_{u_0}:X_{t_0}\to X_{t_0}$ defined by $$T_{u_0}(u)(x,t)=u_0(x)+\int_0^t K_\varphi(u)(x,s)\,\mathrm{d}s\,.$$
Let us show that $T_{u_0}(u)$ remains in $X_{t_0}$, and notice first that boundedness is obvious since $\mu(\R^n)=1$ implies that $\Vert K_\varphi(u)(\cdot,t)\Vert_\infty\leq 2\Vert u(\cdot,t)\Vert_\infty.$

Now, take $x_n\in\Omega\to x\in\bar\Omega$ and $t_n\to t\in[0,t_0]$. Since $u$ is continuous up to the boundary, to pass to the limit and check that $T_{u_0}(u)$ is continuous, we need only to investigate the term 
$$\int_0^{t_n}\int_{\{(x_n+z)\in\Omega\}} u(x_n+z,s)\dmu(z)\mathrm{d}s+\int_0^{t_n}\int_{\{(x_n+z)\in\partial^\mu\Omega\}} \varphi(x_n+z,s)\dmu(z)\mathrm{d}s\,.$$ 
To treat both terms at the same time, consider the extension $\tilde{u}$ defined in \eqref{extension}. Then $\tilde{u}(x_n+z,s)\to \tilde{u}(x+z,s)$ for any $z$ such that $(z+x)\not\in\partial\Omega$ (recall that $u$ and $\varphi$ are continuous on either side of the boundary of $\Omega$), so that $\tilde{u}(x_n+\cdot)\to \tilde{u}(x+\cdot)$ almost everywhere in $\R^n$ and $\tilde{u}$ is uniformly bounded. By dominated convergence we get that $$\int_0^{t_n} \tilde{u}(x_n+z)d\mu(z)\to\int_0^t \tilde{u}(x+z)d\mu(z)\,,$$ which proves that $T_{u_0}(u)$ is continuous in all $[0,t_0]\times\bar\Omega$.

Now it is easy to check that (recall again that $\mu(\R^n)=1$)  $$\sup_{\Omega\times[0,t_0]}|T_{u_0}(v)-T_{u_0}(w)|\le 2t_0\cdot\sup_{\Omega\times[0,t_0]}|v-w|\,,$$ which implies that $T_{u_0}$ is a strict contraction in $X_{t_0}$ provided $t_0<1/2$. Hence in this case, we may use Banach's fixed point theorem in $X_{t_0}$ to obtain a unique solution of $u=T_{u_0}(u)$ in $X_{t_0}$. Once a solution is constructed on $[0,t_0]$, it remains to iterate the procedure and finally get a unique bounded, continuous solution in $\bar\Omega\times[0,\infty[$, which ends the proof.
\end{proof}

\begin{proposition}\label{propboundary}
For any solution $u\in X$ of $(\mathcal{P})$, the equation $u_t=K_\varphi(u)$ holds at any boundary point $(x,t)\in\partial\Omega\times[0,\infty)$, and also at $t=0$ in $\Omega$ (with continuous extension of $\partial_t u$ at $t=0$).
\end{proposition}

\begin{proof}This is done exactly as in the Theorem above : if $(x_n,t_n)\to(x,t)\in\bar\Omega\times[0,\infty)$, then we can pass to the limit in the different terms of the integral equation \eqref{eq0} so that we get the same equation at $(x,t)$, even if it is a boundary point $x\in\partial\Omega$, or $t=0$. 
\end{proof}

\begin{proposition}
Let $\mu\ge0$ be a $L^1(\R^n)$ function and assume that $u\in X$ is a solution of $(\mathcal{P})$ with initial data $u_0$ and boundary data $\varphi$ ; $v\in X$ is a solution with initial data $v_0$ and boundary data $\psi$. Then if $u_0\le v_0$ and $\varphi\le\psi$, we have $u\le v$ in all $\bar\Omega\times[0,\infty[$.
\end{proposition}

\begin{proof}
Let $w_c(x,t)=v(x,t)-u(x,t)+c/(T-t)$ defined in $\bar\Omega\times[0,T]$, for some arbitrary $c,T>0$. The continuous function $w_c$ reaches its min at some point $(x_0,t_0)\in\bar\Omega\times[0,T]$ and obviously $t_0$ cannot be $T$. If $t_0=0$ then we are done, so it remains to investigate the case when $x_0\in\bar\Omega$ and $t>0$. Using the fact that necessarily the equation holds at such a point (even if it is located at the boundary), one gets the following :
$$\partial_t w_c(x_0,t_0)=0=\int_{\R^n}\big\{\tilde{w}_c(x_0+z,t_0)-w_c(x_0,t_0)\big\}\dmu(z)+\frac{c}{(T-t)^2}\,,$$
where $\tilde{w}_c=\tilde{v}-\tilde{u}+c/(T-t)$ (natural extensions of $u$ and $v$ defined by \eqref{extension}). Remember first that $\tilde{w}_c=\psi-\varphi+c(T-t)\ge0$ in $\partial^\mu\Omega\times[0,T]$, and second that $w_c$ reaches its min in $\bar\Omega\times[0,T]$ at $(x_0,t_0)$, so we arrive at 
\begin{eqnarray*}
w_c(x_0,t_0) &= &\int_{\R^n}\tilde{w}_c(x_0+z,t_0)\dmu(z)+\frac{c}{(T-t)^2}\\
&\ge& \int_{\{(x_0+z)\in\Omega\}}w_c(x_0+z,t_0)\dmu(z)\ge +\frac{c}{(T-t)^2}\\
&\ge & \alpha w_c(x_0,t_0)+\frac{c}{(T-t)^2}\,,
\end{eqnarray*}
for some $\alpha\in[0,1[$ (recall that we assume $\mu(\R^n)=1$ and note that the case $\alpha=1$ leads directly to a contradiction since $c>0$). Hence $w_c(x_0,t_0)\ge 0$ and passing to the limit as $c\to0$ yields the result : $u\ge v$ in $\bar\Omega\times[0,T]$. Since $T>0$ is arbitrary, the comparison result holds. 
\end{proof}

\section{The Vanishing Viscosity Approach}

In this section we try another approach to construct a solution, namely using the vanishing viscosity method. It consists in adding a viscous term $-\eps\Delta u$ to the equation to get a classical solution that does take the boundary data, and then pass to the limit as $\eps\to0$. We shall do this in the case $\varphi=0$ for simplicity, the modifications when $\varphi\neq0$ is continuous being rather straightforward.

Recall that in the case $\varphi=0$, we denote by $K_0(u)$ the non-local operator associated with the zero boundary data :
$$K_0(u)(x)=\int_{\{(x+z)\in\Omega\}}u(x+z)\dmu(z)-u(x)\,.$$

\begin{proposition}
Let $\eps>0$. Then for any $u_0\in\mathrm{C}_0(\Omega)$, there exists a unique solution $u=u_\eps$ in $\mathrm{C}(\bar\Omega\times[0,\infty))$ of the following problem :
$$(\mathcal{P}_\eps)\quad\begin{cases}u_t-\eps\Delta u = K_0(u) & \text{in }\Omega	\times(0,\infty)\cr 
u(x,t)=0 & \text{on }\partial\Omega\times(0,\infty)\cr 
u(x,0)=u_0(x) & \text{in }\Omega
\end{cases}$$
\end{proposition}

\begin{proof}
Existence of a solution follows from standard arguments once we have a comparison result and some bounds. In fact the comparison principle for this viscous equation is clear: it is obtained like in the case $\eps=0$ (see previous section) with standard arguments since the viscous term does not change anything here. Moreover, positive constants are solutions (with positive boundary values)  so that for $u_0\in\mathrm{C}_0(\Omega)$ we obtain a solution $u$ which is bounded by $\Vert u_0\Vert_\infty$. Thus $u$ satisfies $u_t-\Delta u =f$ with $f$ bounded, which implies that the zero boundary data is taken in the classical sense, and that $u\in\mathrm{C}(\bar\Omega\times[0,\infty))\cap\mathrm{C}^2(\Omega\times(0,\infty)).$ Uniqueness follows also from the comparison principle. 
\end{proof}

\begin{theorem}
As $\eps\to0$, $u_\eps$ converges pointwise and in $L^\infty(\Omega\times(0,\infty))$-weak\,$\star$ to the unique solution $u\in\mathrm{C}(\bar\Omega\times[0,\infty))$ of the following problem :
$$\quad\begin{cases}
   u_t=K_0(u) & \mbox{in}\ \Omega\times(0,\infty)\,, \cr
   u(x,0)=u_0(x) & \mbox{in}\ \Omega\,.
\end{cases}
$$  
\end{theorem}

\begin{proof}
Let us begin by using a weak formulation for the viscous problem : for any $\phi\in\mathrm{C}^2_0(\Omega\times(0,\infty))$ (\textit{i.e.}, regular and compactly supported in $\Omega\times(0,\infty)$), we have 
$$\int_{\R^n} u_\eps(t)\phi(t)-\int_{\R^n}u_0\phi(0)-\eps\int_0^t\int_{\R^n}u_\eps\Delta\phi\,\mathrm{d}x\mathrm{d}s=\int_0^t\int_{\R^n} K_0(u_\eps)\phi\,\mathrm{d}x\mathrm{d}s\,.
$$
Then since the sequence $(u_\eps)_\eps$ is bounded in $L^\infty$ as $\eps\to0$, we can extract to get a sequence $(u_{\eps_n})$ converging in $L^\infty$-weak\,$\star$ to some function $u_*$, bounded in $\Omega\times(0,\infty)$. This allows us to pass to the limit in the weak fomulation:
$$\int_{\R^n}u_*\phi(t)-\int_{\R^n}u_0\phi(0)=\int_0^t\int_{\R^n} K_0(u_*)\phi\,\mathrm{d}x\mathrm{d}s\,.
$$
This equality proves that at least in the sense of distributions in $\Omega\times(0,\infty)$, the bounded function $u_*$ satisfies $\partial_t u_* = K_0(u_*)$. Now, let us show that $K_0(u_*)$ is a continuous function in $\Omega\times(0,\infty)$, so that the equation holds in the classical sense. Since $u_*$ is bounded and $\mu\in L^1(\R^n)$, we can use the density of continuous functions in $L^1$ to prove that
$$(x,t)\mapsto\int_{\{(x+z)\in\Omega\}}u_*(x+z,t)\dmu(z)\quad\text{is continuous in }\Omega\times(0,\infty)\,.$$ Indeed, this is just using the dominated convergence theorem for a continuous integrand as in the proof of Theorem \ref{thmexun}, and then passing to the limit in $L^1$ by a density argument. Then, since $u_t$ is a continous function in $\Omega\times(0,\infty)$, we can integrate in time the equation satisfied by $v(x,t)=e^tu(x,t)$ to get :
$$v(x,t)-v(x,0)=\int_0^t\int_{\{(x+z)\in\Omega\}}v(x+z,s)\dmu(z)\mathrm{d}s\,,
$$
hence $v$ is continuous and $u_*$ also, in all $\Omega\times(0,\infty)$. Finally, as we already did in the previous section, we can pass to the limit as $x$ goes to the boundary to conclude that the equation is satisfied in $\bar\Omega\times(0,\infty)$. Moreover, $u_*$ can be extended to get a solution $u\in\mathrm{C}(\bar\Omega\times[0,\infty))$, so that $u$ is the unique solution of the problem (by Theorem \ref{thmexun}).

Let us now prove that the limit holds pointwise : take $T>0$ fixed. Since $u_{\eps_n}\to u$ in $L^\infty$-weak\,$\star$, we can always assume that the limit also holds in $L^2$-weak. Let us make some energy estimates in $L^2$ for the solution : since $u_{\eps_n}=0$ on the boundary,
$$\int_\Omega |u_{\eps_n}|^2(t)+\eps_n\int_0^T\int_\Omega|Du_{\eps_n}|^2=\int_\Omega |u_0|^2(t)+\int_0^T\int_\Omega K_0(u_{\eps_n})u_{\eps_n}\,,$$
so that in the limit we get $$\int_\Omega |u_*|^2\le\int_\Omega |u_0|^2+\int_0^T\int_\Omega K_0(u_*)u_*\,.$$ But since this inequality is in fact an equality (obtained by multipling the limit equation by $u_*$ and integrating over $\Omega\times(0,T)$), we get that over all $\Omega\times(0,T)$, $$\lim_{\eps_n\to0}\Vert u_{\eps_n}\Vert_{L^2(\Omega\times(0,T))}=\Vert u_{*}\Vert_{L^2(\Omega\times(0,T))}\,.$$
So, since we have both the weak $L^2$ convergence and the convergence in the $L^2$-norm, it turns out that the convergence holds in $L^2$ strong. Thus, up to another extraction, $u_{\eps_n}\to u_*$ pointwise in $\Omega\times(0,\infty)$ ($T>0$ is arbitrary). Finally, the arguments above show that for any other extraction of the sequence $(u_\eps)$, the limit is necessarily $u_*$ which implies that in fact it is all the sequence $(u_\eps)$ that converges to $u_*$, pointwise in $\Omega\times(0,\infty)$. Thus, $u_\eps\to u$ in $\Omega\times[0,\infty)$ since $u$ is the continuous extension of $u_*$ to $\bar{\Omega}\times[0,\infty)$.
\end{proof}

\section{Regularity Questions}

Let us first mention that obviously if $\mu\in L^1(\R^n)$, any continous solution $u(x,t)$ is $\mathrm{C}^1$ in time, which follows directly from the integral formulation: $u_t$ is continuous in $\bar\Omega\times[0,\infty)$.

Now, as far as the space regularity is concerned, we begin recalling that the equation does not have any regularizing effect when $\mu\in L^1$. More precisely, it was shown in \cite{CCR} in the case $\Omega=\R^n$ (and $\mu$ continuous) that if $u$ is the solution with initial data $u_0$, then $u(t)=e^{-t}u_0 + v(t)$ where $v$ is smooth, so that the regularity of $u_0$ is preserved for all times. See also \cite{Al} for estimates in BUC in $\R^n$ (in a more general context).

We will prove that we can control the space modulus of continuity. If $u\in\mathrm{C}(\Omega\times[0,\infty))$, let us denote by $$\omega(\eta,t)=\sup_{(x,y)\in\Omega^2\atop|x-y|<\eta}|u(x,t)-u(y,t)|\,,$$ and let us begin with a result in the entire space $\R^n$:

\begin{theorem}
Let $\Omega=\R^n$, $\mu\in L^1(\R^n)$ and $u\in\mathrm{C}(\R^n\times[0,\infty))$, a solution with initial data $u_0$ having a modulus of continuity $\omega_0$. Then we have : for any $\eta>0$ and $t>0$, $$\omega(\eta,t)\le\omega_0(\eta)\,.$$
\end{theorem}

\begin{proof}
Let us integrate in time the equation for $v(x,t)=e^tu(x,t)$ and denote by $\omega_v$ the modulus of continuity for $v$ : 
$$v(x,t)-v(y,t) = v(x,0)-v(y,0)+\int_0^t\int_{\R^n}\big\{v(x+z,s)-v(y+z,s)\big\}\dmu(z)\mathrm{d}s\,,
$$
so that if $|x-y|<\eta$, we can estimate the modulus of continuity as follows (recall that by assumption, $\mu(\R^n)=1$) :
$$|v(x,t)-v(y,t)|\le\omega_0(\eta)+\int_0^t\omega_v(\eta,s)\mathrm{d}s\,.$$
Hence using Gronwall's lemma we get that $\omega_v(\eta,t)\le\omega_0(\eta)e^t$, and going back to $u(x,t)=e^{-t}v(x,t)$ gives the result.
\end{proof}

In this generality, the Theorem is optimal : take $N=1$ and for instance $\mu(z)=1/2$ on $[-1;1]$, $\mu(z)=0$ if $|z|>1$. Then if $u_0(x)=ax$, the solution is stationnary and so the modulus of continuity is constant. 

The situation is a little bit more complicated in a bounded domain $\Omega$ since the boundary terms may induce some sort of discontinuity. So let us first restrict ourselves to the case $\varphi=0$ and let us introduce some notation :

For any $x\in\R^n $, let $\tau_x^{-1}\Omega=\{z:(x+z)\in\Omega\}$ and for any two sets $A,B$, let $A\,\Delta\,B = (A\setminus B)\cup(B\setminus A)$. Now, for any $\eta>0$, let us define the following quantities : 
\begin{eqnarray*} & \lambda(\eta)=\sup\limits_{(x,y)\in\Omega^2\atop|x-y|<\eta}\mu\big(\tau_{x}^{-1}\Omega\cap\tau_{y}^{-1}\Omega\big)\,,\\
& \gamma(\eta)=\sup\limits_{(x,y)\in\Omega^2\atop|x-y|<\eta}\mu\big(\tau_x^{-1}\Omega\,\Delta\,\tau_y^{-1}\Omega\big)\,.
\end{eqnarray*}
Notice that obviously, $\lambda(\eta),\gamma(\eta)\in[0,1]$ and that $\gamma(\eta)$ goes to zero as $\eta\to0$. Roughly speaking, $\lambda(\eta)$ is related to the total mass of $\mu$ in $\widetilde{\Omega}=\bigcup\limits_{x\in\Omega}\tau_x^{-1}\Omega$, while $\gamma(\eta)$ is related to the mass of $\mu$ involved at the boundary. 

\begin{theorem}
Let us assume that $\Omega$ is bounded and regular, $\mu\ge0$ is a $L^1$ measure, $\varphi=0$ and that $u$ is a solution with initial data $u_0\in\mathrm{C}_0(\bar\Omega)$. Then the following estimates hold : 

\noindent\textbullet\; If $\lambda(\eta)<1$, then for any $t>0$
$$\omega(\eta,t)\le\Big(\omega_0(\eta)+\gamma(\eta)\cdot\Vert u_0\Vert_\infty\frac{e^{(1-\lambda(\eta)) t}-1}{1-\lambda(\eta)}\Big)e^{(\lambda(\eta)-1)t}\,.$$
\noindent\textbullet\; If $\lambda(\eta)=1$, then for any $t>0$, 
$$\omega(\eta,t)\le\Big(\omega_0(\eta)+\gamma(\eta)\cdot\Vert u_0\Vert_\infty\cdot t\Big)\,.$$
\end{theorem}

\begin{proof}
The proof proceeds with the same ideas as for the case $\Omega=\R^n$ : integrating in time the equation for $v(x,t)=e^tu(x,t)$ and again denoting by $\omega_v$ the modulus of continuity for $v$ one gets : 
\begin{eqnarray*}v(x,t)-v(y,t) & = & \int_0^t\int_{\tau_x^{-1}\Omega}v(x+z,s)\dmu(z)\mathrm{d}s-\int_0^t\int_{\tau_y^{-1}\Omega} v(y+z,s)\dmu(z)\mathrm{d}s\,,
\end{eqnarray*}
so that if $|x-y|<\eta$, we can estimate the modulus of continuity as follows :
$$|v(x,t)-v(y,t)|\le\int_0^t\int_{\tau_x^{-1}\Omega\cap\tau_y^{-1}\Omega}\big|v(x+z,s)-v(y+z,s)\big|\dmu(z)\mathrm{d}s + R(x,y,t)\,$$
where $R(x,y,t)$ contains the integral over $\big(\tau_x^{-1}\Omega\,\Delta\,\tau_y^{-1}\Omega\big)$, that we estimate uniformly using the fact that $\Vert v(s)\Vert_\infty\le\Vert u_0\Vert_\infty\cdot e^s$ : $$R(x,y,t)\le \Vert u_0\Vert_\infty \cdot\gamma(\eta)\cdot\int_0^t e^s\,\mathrm{d}s\,.$$ 
Thus, going back to the estimate, we have for any $|x-y|<\eta$, 
$$|v(x,t)-v(y,t)|\le\omega_v(\eta,0)+\int_0^t\big\{\lambda(\eta)\cdot\omega_v(\eta,s)+\gamma(\eta)\cdot\Vert u_0\Vert_\infty\cdot e^s \big\}\,\mathrm{d}s \,.
$$
Passing to the supremum over all $(x,y)$ such that $|x-y|<\eta$, and using Gronwall's lemma we obtain: 
$$\omega_v(\eta,t)\le \big(\omega_v(\eta,0) + \gamma(\eta)\cdot\Vert u_0\Vert_\infty\frac{e^{(1-\lambda(\eta)) t}-1}{1-\lambda(\eta)}\big)e^{\lambda(\eta)t}\,, 
$$
and finally going back to $u(x,t)=e^{-t}v(x,t)$ yields the desired estimate. When $\lambda(\eta)=1$, then the estimate is easier to obtain, using also Gronwall's lemma.
\end{proof}

Before stating a corollary, let us give a heuristic idea of the meaning of those estimates: if $\lambda=1$, we see that the estimate on the modulus of continuity deteriorates as time increases, because of the discontinuity at the boundary (the jump between $u$ and $\varphi=0$) which propagates in the interior. On the contrary, if $\lambda(\eta)<1$, this effect is compensated by the fact that the measure $\mu$ is not fully integrated, so that we gain an exponential decay of the modulus of continuity : some kind of regularization occurs ! (See end Section for some examples).

Let us give two sufficient conditions that explain when both cases may appear (recall that $\widetilde{\Omega}=\bigcup\limits_{x\in\Omega}\tau_x^{-1}\Omega$) :

\begin{proposition}
If $\mu(\R^n\setminus\widetilde{\Omega})>0$, then for any $\eta>0$, $$\lambda(\eta)\le 1-\mu(\R^n\setminus\widetilde{\Omega})<1\,.$$ On the other hand, if there exists $x\in\Omega$ such that $\mathrm{supp}(\mu)\subset\tau_x^{-1}\Omega$ (with strict inclusion), then $\lambda(\eta)=1$ for $\eta>0$ sufficiently small. 
\end{proposition}

\begin{proof}
Let us first notice that obviously, for any $x,y\in\Omega$, $\tau_x^{-1}\Omega\cap\tau_y^{-1}\Omega\subset\widetilde{\Omega}$ so that $\mu(\tau_x^{-1}\Omega\cap\tau_y^{-1}\Omega)\le\mu(\widetilde{\Omega})$, which is independent of $\eta>0$ and $x,y\in\Omega$. Thus passing to the supremum and using the fact that the total mass of $\mu$ is 1, one gets the first result.

Now for the second statment: if such an $x$ exists, then for $\eta$ small enough, any $y\in\Omega$ such that $|x-y|<\eta$, will satisfy : $\mathrm{supp}(\mu)\subset\tau_x^{-1}\Omega\cap\tau_y^{-1}\Omega$. Hence $\mu(\tau_x^{-1}\Omega\cap\tau_y^{-1}\Omega)=1$ and passing to the supremum for such $x,y$ yields the result.
\end{proof}

We refer to the final Section where explicit examples are given, by let us mention that more or less, $\lambda>1$ if $\mu$ has a huge support with respect to the shape of $\Omega$, whereas $\lambda=1$ if $\mu$ is rather concentrated (again with respect to the shape of $\Omega$). The limit being when $\mu=\delta_0$ is so concentrated that nothing moves. 

\begin{remark}
If $\varphi\neq0$ is bounded and continuous, the estimate involves another term depending on the integral of $\varphi$ on $A_\mu(\eta)=\partial^\mu\Omega_\eta\setminus\partial^\mu\Omega$ : if 
$$\theta(\eta,t)=\max_{\Omega\times(0,t)}\int_{\{(x+z)\in A_\mu(\eta)\}}\varphi(x+z)\dmu(z),
$$ then we can estimate the modulus of continuity using again Gronwall's Lemma, taking into account $\theta(\eta,t)$. The details are left to the reader.
\end{remark}

\section{Qualitative Properties of Solutions}

\begin{theorem}\label{thmspeed}
Let us assume that $\mu\ge0$ is a measure such that for some $\eta>0$, $$\mathrm{supp}\,(\mu)\supset\{|z|\le\eta\}\,.$$ 
Then if $u(x,t)\in X$ is a solution of problem $(\mathcal{P})$ with $\varphi=0$ and $u_0\in\mathrm{C}(\Omega)$ non-negative and compactly supported in $\Omega$, we have $$\forall\,(x,t)\in\bar\Omega\times(0,\infty),\quad u(x,t)>0\,.$$
\end{theorem}

\begin{proof}
Let us first observe that if we set $v(x,t)=e^t\tilde{u}(x,t)$ where $\tilde{u}$ is the extension of $u$ defined by \eqref{extension}, then $v$ satisfies $$\partial_t v=\int_{\R^n}v(x+z)\dmu(z)\ge0\,,$$ so that once $v$ (hence $u$ also) is positive, it stays as such for later times. Let $A=\mathrm{supp}\,(u_0)\subset\Omega$ which is compact and notice first that obviously $u\ge0$ in $\Omega\times[0,\infty)$, using the comparison principle. Now, let $x_0\in\bar\Omega\setminus A$, and assume that $\mathrm{dist}\,(x_0;A)<\eta$. Using the equation at $(x_0,0)$ (remember that the equation holds even on the boundary), one gets
$$v_t(x_0,0)\ge\int_{|z|<\eta}v_0(x+z,0)\dmu(z)>0\,.$$
Thus we get that necessarily, for any $t>0$, $u(x,t)$ will be positive in 
$A_\eta=\{x\in\bar\Omega:\mathrm{dist}\,(x,A)<\eta\}$. Iterating the same argument at $t>0$ arbitrary small, we get that $u(x,t)>0$ in $A_{k\eta}$ for any $k\in\mathbb{N}$. Thus it is clear that $u(x,t)>0$ for any $x\in\bar\Omega$, for any time $t>0$.
\end{proof}

As a corollary of this result, we have obtained that the boundary data is never taken, more precisely :

\begin{corollary}
In the vanishing viscosity method, a boundary layer occurs and the limit solution is strictly positive on $\partial\Omega$ for any time $t>0$. Thus the boundary data is never taken in general.
\end{corollary}

\section{Conclusions, Remarks and Examples}

We have proposed a general framework to deal with Dirichlet problems for the non-local equation \eqref{eq0}. If $\mu\in L^1(\R^n)$, we have proved existence and uniqueness of a continuous solution for any continuous initial and boundary data. Furthermore, under some condition on the support of $\mu$ near the origin, we have also proved a positivity property, and that the boundary values will not be taken in the classical sense. Let us end this note with some remarks and extensions :

\

\textbf{1.} If $\mu\in L^\infty$, a similar existence and uniqueness result as in Theorem \ref{thmexun} may be obtained for data $u_0,\varphi$ in $L^1$, using similar fixed-point arguments. But if $\mu$ is a finite measure with singular part (like Dirac deltas), then the above construction does not work since the integral term may not be continuous. However, we think that similar things could maybe be done for measures less concentrated than the delta. 

\textbf{2.} An interesting feature of the equation is that it holds true everywhere on $\partial\Omega$ (see Proposition \ref{propboundary}). This explains why, although the boundary data is not taken in the usual sense, there would not be any problem with the viscosity definition of boundary values. Indeed, passing to the limit in the viscous equation can be done also in the framework of viscosity solutions since the convergence holds in $L^\infty$-weak$\star$.

\textbf{3.} Theorem \ref{thmspeed} shows that propagation with infinite speed occurs like for the Heat Equation, although the phenomenon is of different nature (due to the non-local term and not instantaneous diffusion as it is the case for the heat equation). But even more, $u$ also becomes positive at the boundary $x\in\partial\Omega$. So, even if both $\varphi$ and $u_0$ are continuous, in general $u(x,t)$ will not take the boundary data on $\partial\Omega$ in the usual sense (this was proved for stationary solutions in \cite{CCR}). And again, we recover the necessity to define a suitable notion of boundary for non-local operators, $\partial^\mu\Omega$ here. 

\textbf{4.} Note that some property on the support of $\mu$ near the origin is necessary to have this propagation property, although our criteria is maybe not optimal. Indeed, if $\mu$ is too much off-centered, then the operator will not see anything in $\Omega$ itself: consider for instance $\Omega=]-10;+10[$, $\varphi=0$ and let $\mu$ be the caracteristic function of the set $]-1,0[$. If $u(0)\in\mathrm{C}_0(\Omega)$, it will remain compactly supported on the left for any $t>0$.

\textbf{5.} We showed a kind of regularizing effect in the case of bounded domains, provided the measure $\mu$ satfisfies $\lambda(\eta)<1$. Let us consider the case $\Omega=]-1;1[$ and let us give an example of a measure $\mu$ satisfying this property : $\mu(z)=ce^{-|z|^2}$, where $c>0$ is such that $\mu(\R)=1$. Then clearly, since $\widetilde{\Omega}=]-2;2[$, $\lambda(\eta)\le 1-\int_{-2}^{+2}\dmu(z)<1$. One can also have a good estimate of the rate in the exponential decay. Now if $\mu(z)$ is the characteristic function of the set $]-1/2;1/2[$, we have clearly that for any $\eta<1/2$, $\lambda(\eta)=1$.

\textbf{6.} A different story happens when $\mu$ is not integrable near the origin, which is the case for the fractional Laplace operator. In this case, the boundary data are taken since the measure is so singular that discontinuities at the boundary are not allowed (see  \cite{BCI}).
\

\

{\small
\section*{Appendix : The Probabilistic Interpretation}

In this Appendix, we briefly explain the probabilistic context in which the Cauchy problem may be understood. Let $(\Omega,\mathcal{F},\P)$ be a measured space and $Y_0:\Omega\to\R$ be a random variable with law $u_0\in\ L^1(\R^N)$ (we assume that $\int u_0=1$). Now let $(Y_i)_{i\in\N}$ be a sequence of independent, identically distributed random variables of law $\mu$, and $t\mapsto N_t$ be a random variable on $\N$ with Poisson law of intensity $t>0$ : $$\P(N_t=k)=e^{-t}\frac{t^k}{k!}\,.$$

Then we consider the stochastic process $X_t=Y_0+\sum\limits_{i=1}^{N_t}Y_i$. The term $Y_0$ represents the initial state of the stochastic process $X_t$, while the sum represents a random number of jumps $N_t$, each jump being independent of the others, with probability distribution $\mu$. Such a process is named a \textit{composed Poisson process}, which is a special case of the more general L\'evy processes. If we consider the density $u(t)$ of the probability distribution of $X_t$, it can be easily shown that precisely, $$\frac{\partial u}{\partial t}=\int_{\R^N}\{u(x+z,t)-u(x,t)\}\dmu(z)\,.$$

Now, by definition the characteristic function of $X_t$ is nothing but the Fourier transform of $u(t)$, and it is known that for any L\'evy process, there exists a uniquely determined triplet $(A,\gamma,\nu)$ such that :
$$\mathds{E}(e^{i\xi X_t})=e^{t\varphi(z)}\ \text{with}\ \varphi(z)=-\frac{1}{2}z\cdot Az+\gamma\cdot z + \int_{\R^N}\big(e^{iz\xi}-1-z\xi\mathds{1}_{|z|\leqslant 1}(z)\big)\,\mathrm{d}\nu(z)\,.$$

In this decomposition, $A$ is a positive-definite matrix, $\gamma\in\R^n$ and $\nu$ is the L\'evy measure, satisfying the following properties : $\int_{|z|\leqslant1}z^2\mathrm{d}\nu\,,\ \int_{|z|\geqslant1}\mathrm{d}\nu <\infty$. The infinitesimal generator of the Feller semi-group associated to $(X_t)$ is the following :
$$Lf:=-\frac{1}{2}\sum_{i,j=1}^N a_{ij} \frac{\partial^2 f}{\partial x_i\partial x_j} + \sum_{i=1}^N\gamma_i\frac{\partial f}{\partial x_i}  + \int_{\R^N}\big(f(x+z)-f(x)-\nabla f(x)\cdot z\mathds{1}_{|z|\leqslant1}(z)\big)\,\mathrm{d}\nu(z)\,.$$

In the special case of composed Poisson processes, the matrix $A$ is nul, $\gamma=0$ and the corrector term involving the gradient of $f$ vanishes for symmetric measures.}

\

\end{document}